\documentclass[10pt,twoside,english]{amsart}
\advance\oddsidemargin by -1.0cm
\advance\evensidemargin by -1.0cm
\advance\topmargin by -1.0cm 
\usepackage{amssymb}
\usepackage{babel}
\usepackage{amstext}
\usepackage{amscd}   
\usepackage{epsfig} 
 \usepackage{rotating}
\theoremstyle{plain}

\theoremstyle{definition}

\newcommand{\Ric}{\ensuremath{\mathrm{Ric}}}

\newtheorem{proposition}{Proposition}[section]
\newtheorem{corollary}[proposition]{Corollary}
\newtheorem{lemma}[proposition]{Lemma}

\newtheorem{remark}[proposition]{Remark}
\newtheorem{theorem}[proposition]{Theorem}

\newcommand{\be}{\begin{equation}}
\newcommand{\ee}{\end{equation}}
\newcommand{\bea}{\begin{eqnarray}}
\newcommand{\eea}{\end{eqnarray}}
\newcommand{\bean}{\begin{eqnarray*}}
\newcommand{\eean}{\end{eqnarray*}}

\begin{document}

\title{Some integrability conditions for almost K\"ahler manifolds}
\author{K.-D. Kirchberg}
\date{}
\maketitle \vspace{1cm}

\begin{abstract}
\noindent Among other results, a compact almost
K\"ahler manifold is proved to be K\"ahler if the Ricci tensor is semi-negative
and its length coincides with that of the star Ricci tensor
or if the Ricci tensor is semi-positive and its first order covariant
derivatives are Hermitian. Moreover, it is shown that there are no compact 
almost K\"ahler manifolds with harmonic Weyl tensor and non-parallel
semi-positive Ricci tensor. Stronger results are obtained in dimension 4. \\

\end{abstract}
\bigskip

\setcounter{section}{-1}
\noindent\section{\bf Introduction}

To find suitable curvature conditions that imply the integrability
of the almost complex structure is one of the most important problems
concerning almost K\"ahler manifolds. In this context, the starting 
point for many investigations was Goldberg's conjecture of 1969, which states
that every compact Einstein almost K\"ahler manifold is necessarily K\"ahler 
\cite{14}. Important progress was made by K. Sekigawa in 1987. He proved the
Goldberg conjecture for non-negative scalar curvature \cite{27}. In case of
negative scalar curvature, no proof is known so far. There are attempts to
construct counterexamples against this part of the Goldberg conjecture. Our 
paper deals  with several kinds of curvature conditions that force
an almost K\"ahler manifold to be  K\"ahler, i.e., that the almost complex 
structure of an almost 
K\"ahler manifold is integrable. One of our main results is a
generalization of Sekigawa's theorem mentioned above. We prove that a compact
almost K\"ahler manifold is K\"ahler if the Ricci tensor is semi-positive
$(\Ric \ge 0)$ and its first order covariant derivatives commute with the 
almost complex structure (Corollary \ref{cor-3-2}). In dimension 4, 
the supposition that $\Ric \ge 0$ can be replaced by the weaker condition that 
the star scalar curvature $S_{\star}$ is non-negative (Corollary 
\ref{cor-3-3}). This result is more general than the theorem that  H. Satoh 
\cite{25} proved recently. Satoh's theorem states that every compact almost 
K\"ahler manifold with semi-positive Ricci tensor and harmonic Weyl tensor 
$(\delta W=0)$ is already K\"ahler. We show that there are no compact almost 
K\"ahler manifolds with harmonic Weyl tensor and semi-positive, 
non-parallel Ricci tensor (Corollaries \ref{cor-3-4}, \ref{cor-3-5}).
Thus, the suppositions of Satoh's theorem imply that the Ricci tensor is 
parallel.\\
It is well known that an almost K\"ahler manifold is K\"ahler if its
star scalar curvature coincides with the scalar curvature $S$. A similar
result is our Proposition \ref{prop-3-1}, which states that an almost 
K\"ahler manifold with semi-positive Ricci tensor is already K\"ahler if
the Hermitian parts of the Ricci tensor and the star Ricci tensor have
the same length. In dimension 4, the supposition that $\Ric \ge 0$ can replaced
by the essential weaker condition  that $S \ge 0$ or 
that the set of zeros of $S + S_{\star}$ is
nowhere dense (Proposition \ref{prop-3-2}). For a compact almost K\"ahler
manifold with semi-negative Ricci tensor $(\Ric \le 0)$, the K\"ahler
property is forced by the supposition that the length of the star Ricci
tensor $\Ric_{\star}$ coincides with that of the Ricci tensor (Theorem
\ref{thm-3-1}).\\

In the case of a compact Einstein almost K\"ahler $n$-manifold, the curvature
inequality\\

$\displaystyle (*) \hfill | \tilde{R}^- |^2 + | \Ric_{\star} |^2 \ge
\frac{1}{n} S \cdot S_{\star}$ \hfill \mbox{}\\

forces the K\"ahler property (Theorem \ref{thm-3-2}). Here $\tilde{R}^-$ is a
part of the curvature tensor depending  on the almost complex structure
$J$ and the Weyl tensor $W$ only. This inequality is satisfied trivially
if $S_{\star} \ge 0$. Since $S \ge 0$ implies $S_{\star} \ge 0$, we obtain
Sekigawa's result. With regard to the Goldberg conjecture it may be 
interesting to investigate for which compact Einstein almost K\"ahler
manifolds $(*)$ is valid if $S <0$.\\
In order to obtain the results for the compact case, 
we modify a well known basic
Weitzenb\"ock formula. So we find two integral formulas of different kind 
(Proposition \ref{prop-2-1}). The second one is applicable if a certain
number $Q(J)$ vanishes. $Q(J)$ is a globally defined obstruction against the 
integrability of the almost complex structure $J$ of every compact almost
K\"ahler manifold. We prove that a compact almost K\"ahler manifold with 
semi-positive Ricci tensor is K\"ahler if and only if $Q(J)=0$ (Theorem
\ref{thm-3-3}). This theorem and the corresponding $4$-dimensional version
(Theorem \ref{thm-3-4}) are essential results of this paper. In these
theorems the suppositions that 
$\Ric \ge 0$ and $S_{\star} \ge 0$, respectively, 
can be replaced by weaker curvature inequalities (Remark \ref{remark-3-1}).
Moreover, we list some geometrical conditions, each of which implies $Q(J)=0$ 
(Remark \ref{remark-3-2}).

\section{\bf Preliminaries}

Let $(M, g,J)$ be an almost Hermitian manifold of dimension $n =2m$ with 
Riemannian metric $g$ and almost complex structure $J$. Then, by definition
\begin{equation} \label{gl-1}
J^2 = -1
\end{equation}

and $g$ is $J$-invariant, i.e., we have
\begin{equation} \label{gl-2}
g(JX, JY)=g (X,Y)
\end{equation}

for all vector fields $X, Y$. The corresponding fundamental $2$-form 
$\Omega$ is defined by $\Omega (X,Y):= g(JX,Y)$. An almost Hermitian
manifold is called almost K\"ahler  if its fundamental form is closed
\begin{equation} \label{gl-3}
d \Omega =0 \ . 
\end{equation}

It is well known that the basic equations (1)-(3) of an almost K\"ahler 
manifold imply that $\Omega$ is also co-closed
\begin{equation} \label{gl-4}
\delta \Omega =0
\end{equation}

and that $J$ satisfies the so-called quasi K\"ahler condition
\begin{equation} \label{gl-5}
\nabla_{JX} J= \nabla_X J \circ J , 
\end{equation}

where $\nabla$, as usually, denotes the Levi-Civita covariant derivative
corresponding to $g$. (3) is equivalent to 
\begin{equation} \label{gl-6}
g(( \nabla_X J)Y,Z)+g((\nabla_Y J)Z, X) +g((\nabla_Z J)X, Y)=0 . 
\end{equation}

By (1) and (2), it holds that
\begin{equation} \label{gl-7}
g(JX,Y)= - g (X, JY) , 
\end{equation}

i.e., $J$ is anti-selfadjoint (skew symmetric)
\begin{equation} \label{gl-8}
J^* = - J . 
\end{equation}

Applying $\nabla_X$ to equation (1) we obtain
\begin{equation} \label{gl-9}
\nabla_X J \circ J + J \circ \nabla_X J=0  . 
\end{equation}

In the following we use the notation
\begin{equation*}
\nabla^2_{X,Y} := \nabla_X \circ \nabla_Y - \nabla_{\nabla_X Y}
\end{equation*}

for the tensorial covariant derivatives of second order. Then the Riemannian
curvature tensor $R$ of the metric $g$ is given by
\begin{equation} \label{gl-10}
R(X,Y) Z= \nabla^2_{X,Y} Z- \nabla^2_{Y,X} Z . 
\end{equation}

Moreover, using (5) we obtain the equations
\begin{eqnarray} \label{gl-11}
\nabla^2_{X,JY} J &=& \nabla^2_{X,Y} J \circ J + \nabla_Y J \circ
\nabla_X J - \nabla_{(\nabla_X J) Y } J , \\ \label{gl-12}
\nabla^2_{X,JY} J &=& - J \circ \nabla^2_{X,Y} J - \nabla_X J \circ \nabla_Y J
- \nabla_{(\nabla_X J) Y} J . 
\end{eqnarray}

For endomorphisms $A,B$ of the tangent bundle $TM$, we use the notations
\begin{displaymath}
[A,B] := A \circ B - B \circ A \quad , \quad \{ A,B \} := A \circ B + B \circ A
\end{displaymath}

for their commutator and anti-commutator, respectively. Then, for any 
endomorphism $A$ and the almost complex structure $J$, we have the relations
\begin{equation} \label{gl-13}
[ \{ A , J \} , J]=0 \quad , \quad \{ [ A, J] , J\} =0 . 
\end{equation}

By (11) and (12), we immediately obtain
\begin{equation} \label{gl-14}
\{ \nabla^2_{X,Y} J, J \} = - \{ \nabla_X J, \nabla_Y J \} . 
\end{equation}

Let $(X_1 , \ldots , X_n)$ be any local frame of vector fields on $M$. Then, by
$(X^1 , \ldots , X^n)$ we denote the associated coframe, which, using the
convention of summation, is defined by $X^k := g^{kl} X_l$, where $(g^{kl})$
is the inverse of the matrix $(g_{kl})$ with $g_{kl} := g(X_k , X_l)$. Thus,
in the 
case of an orthonormal frame, we have $X^k = X_k (k=1, \ldots , n)$. In the
following we sometimes use orthonormal frames. We remark that (4) is then
locally equivalent to
\begin{equation} \label{gl-15}
(\nabla_{X_k} J) X^k =0 , 
\end{equation}

implying the equation
\begin{equation} \label{gl-16}
(\nabla^2_{X, X_k} J) X^k =0
\end{equation}

for any vector field $X$. The Ricci tensor is given by
\begin{equation} \label{gl-17}
\Ric (X) := R(X, X_k) X^k
\end{equation}

and the star Ricci tensor of the almost K\"ahler manifold $(M, g,J)$ is
defined by
\begin{equation} \label{gl-18}
\Ric_{\star} (X) := R(JX, JX_k)X^k . 
\end{equation}

Moreover, we use the notations
\begin{eqnarray}
\Ric^+ & := & \frac{1}{2} (\Ric - J \circ \Ric \circ J) = - \frac{1}{2} J \circ
\{ \Ric , J \} , \label{gl-19} \\
\Ric^- & := & \frac{1}{2} (\Ric + J \circ \Ric \circ J)= \frac{1}{2} J \circ
[ \Ric , J ] , \\
\Ric^+_{\star} & := & \frac{1}{2} (\Ric_{\star} - J \circ \Ric_{\star} \circ J
)= - \frac{1}{2} J \circ \{ \Ric_{\star} , J \} , \label{gl-21}\\
\Ric^-_{\star} & := & \frac{1}{2} (\Ric_{\star} + J \circ \Ric_{\star} \circ J)
= \frac{1}{2} J \circ [ \Ric_{\star} , J ] . 
\end{eqnarray}

By definition, we have
\begin{equation} \label{gl-23}
\Ric = \Ric^+  + \Ric^- \quad , \quad \Ric_{\star} = \Ric^+_{\star}
+ \Ric^-_{\star}
\end{equation}

and from (13) and (19) - (22) we see that
\begin{equation} \label{gl-24}
[ \Ric^+ , J]= \{ \Ric^- , J \} = [ \Ric^+_{\star} , J]= \{ \Ric^-_{\star} , 
J \} =0 . 
\end{equation}

Obviously, the endomorphisms $\Ric^+$ and $\Ric^-$ are symmetric
\begin{equation} \label{gl-25}
(\Ric^{\pm} )^* = \Ric^{\pm} . 
\end{equation}

Using the first Bianchi identity we find
\begin{equation} \label{gl-26}
\Ric_{\star} = \frac{1}{2} R(X_k , JX^k) \circ J , 
\end{equation}

which implies
\begin{equation} \label{gl-27}
(\Ric_{\star} )^* = - J \circ \Ric_{\star} \circ J . 
\end{equation}

From (21), (22) and (\ref{gl-27}) we see that $\Ric^+_{\star}$ is symmetric
and $\Ric^-_{\star}$ is skew symmetric
\begin{equation} \label{gl-28}
(\Ric^{\pm} )^* = \pm \Ric^{\pm}_{\star} . 
\end{equation}

We remark that in the K\"ahler case $(\nabla J=0)$ we have $\Ric_{\star} = 
\Ric$ and $\Ric^-_{\star} = \Ric^- =0$. The Ricci form $\rho$ and the star
Ricci form $\rho_{\star}$ are defined by
\begin{displaymath}
\rho (X,Y):=g((J \circ \Ric^+)X, Y) \quad \mbox{and} \quad \rho_{\star}
(X,Y):= g(( J \circ \Ric^+_{\star} ) X, Y ) , 
\end{displaymath}

respectively. Both Ricci forms are Hermitian ($J$-invariant)
\begin{equation} \label{gl-29}
\rho (JX , JY)= \rho (X,Y) \quad , \quad \rho_{\star} (JX , JY) = 
\rho_{\star} (X,Y) . 
\end{equation}

Besides the scalar curvature $S:= \mathrm{tr} (\Ric) = \mathrm{tr} (\Ric^+)$
also the star scalar curvature $S_{\star} := \mathrm{tr} (\Ric_{\star})=
\mathrm{tr} (\Ric^+_{\star})$ is considered. Further, for all $X,Y \in \Gamma
 (TM)$, we have the curvature endomorphism $\tilde{R} (X,Y)$ defined by
\begin{displaymath}
\tilde{R} (X,Y) := \frac{1}{4} [R (X,Y) - R(JX, JY) , J] \circ J . 
\end{displaymath}

We see that $\tilde{R}$ has  the properties
\begin{equation} \label{gl-30}
\tilde{R} (X,Y,)^* =- \tilde{R} (X,Y) = \tilde{R} (Y, X) , 
\end{equation}
\begin{equation} \label{gl-31}
\{ \tilde{R} (X,Y) , J \} = 0 , 
\end{equation}
\begin{equation} \label{gl-32}
\tilde{R} (JX, JY)= - \tilde{R} (X,Y) . 
\end{equation}

$\tilde{R} (X,Y)$ decomposes as follows
\begin{equation} \label{gl-33}
\tilde{R} (X,Y) = \tilde{R}^+ (X,Y) + \tilde{R}^- (X,Y) , 
\end{equation}

where $\tilde{R}^+$ and $\tilde{R}^-$ are given by
\begin{displaymath}
\tilde{R}^{\pm} (X,Y) := \frac{1}{2} (\tilde{R} (X,Y) \pm 
\tilde{R} (JX, Y) \circ J) . 
\end{displaymath}

Obviously, $\tilde{R}^+$ and $\tilde{R}^-$ also have the properties
(\ref{gl-30}) - (\ref{gl-32}) and it holds that
\begin{equation} \label{gl-34}
\tilde{R}^{\pm} (JX, Y) \circ J = \pm \tilde{R}^{\pm} (X,Y) . 
\end{equation}

Using (11), (12) a straightforward calculation yields the relation
\begin{equation} \label{gl-35}
\tilde{R}^+ (X,Y) = \frac{1}{4} \nabla_{\varphi (X,Y)} J , 
\end{equation}

where $\varphi$ is the vector-valued 2-form defined by $\varphi (X,Y):=
(\nabla_X J)Y - (\nabla_Y J)X$. Using (5) and (9) it is easy to see
that $\varphi$ has the properties
\begin{equation} \label{gl-36}
\varphi (JX, JY)= - \varphi (X,Y) , 
\end{equation}
\begin{equation} \label{gl-37}
\varphi (JX, Y)= \varphi (X, JY) = - J(\varphi (X, Y)) . 
\end{equation}

Furthermore, from (6) we derive
\begin{equation} \label{gl-38}
\nabla_{\varphi (X,Y)} J= - g((\nabla_{X_k} J)X,Y) \cdot \nabla_{X^k} J . 
\end{equation}

Inserting (\ref{gl-38}) into (\ref{gl-35}) we obtain Gray's identity 
\cite{15}
\begin{equation} \label{gl-39}
\tilde{R}^+ = - \frac{1}{4} \nabla_{X_k} \Omega \otimes \nabla_{X^k} J . 
\end{equation}

It is well known that $\tilde{R}^-$ is already determined by the Weyl tensor.
Let $\widetilde{\Ric}$ be the endomorphism of $TM$ locally defined by
$\widetilde{\Ric} (X):= \tilde{R} (X, X_k)X^k$. Using (\ref{gl-33}) and
(\ref{gl-34}) we obtain
\begin{equation} \label{gl-40}
\widetilde{\Ric} (X)= \tilde{R}^+ (X, X_k) X^k . 
\end{equation}

By (\ref{gl-39}) and (\ref{gl-40}), we find
\begin{equation} \label{gl-41} 
\widetilde{\Ric} = - \frac{1}{4} \nabla_{X_k} J \circ \nabla_{X^k}J . 
\end{equation}

Thus, $\widetilde{\Ric}$ is symmetric and semi-positive
\begin{equation} \label{gl-42}
(\widetilde{\Ric})^* = \widetilde{\Ric} , 
\end{equation}
\begin{equation} \label{gl-43}
\widetilde{\Ric} \ge 0 . 
\end{equation}

Moreover, $\widetilde{\Ric}$ commutes with the almost complex structure
\begin{equation} \label{gl-44}
[ \widetilde{\Ric} , J]=0 . 
\end{equation}

In the following we use the Bochner Laplacian $\nabla^* \nabla$ locally
given by $\nabla^* \nabla := - \nabla^2_{X_k , X^k}$. From (14) and
(\ref{gl-42}) we immediately obtain
\begin{equation} \label{gl-45}
\{ \nabla^* \nabla J, J \} = - 8 \widetilde{\Ric}
\end{equation}

and using (\ref{gl-16}) we find
\begin{equation} \label{gl-46}
(\nabla^2_{X_k , X} J)X^k =( J \circ \Ric - \Ric_{\star} \circ J) X . 
\end{equation}

Further, (\ref{gl-6}) implies the equation
\begin{equation} \label{gl-47}
g((\nabla^2_{V,X} J)Y,Z)+g((\nabla^2_{V,Z} J)X,Y)+g((\nabla^2_{V,Y} J)Z,X)=0 , 
\end{equation}

which is valid for all vector fields $V, X, Y, Z$. Contracting this
equation and using (\ref{gl-19}), (\ref{gl-27}), (\ref{gl-46}) we obtain
\begin{equation} \label{gl-48}
\nabla^* \nabla J = 2 (\Ric_{\star} - \Ric^+) \circ J . 
\end{equation}

By (\ref{gl-45}), (\ref{gl-48}) and (\ref{gl-21}), we find the identity
\begin{equation} \label{gl-49}
\widetilde{\Ric} = \frac{1}{2} (\Ric^+_{\star} - \Ric^+) . 
\end{equation}
 
Multiplying (\ref{gl-46}) by $J$ and using (23), (\ref{gl-27}) and (
\ref{gl-28}) we obtain
\begin{equation} \label{gl-50}
(J \circ \nabla^2_{X_k, X} J) X^k =(\Ric^+_{\star} - \Ric^+)X-
(\Ric^-_{\star} + \Ric^-)X . 
\end{equation}

For endomorphisms $A,B$ of the tangent bundle $TM$, we use the scalar
product $\langle A, B \rangle := \mathrm{tr} (A \circ B^*)$. On the other 
hand, the scalar product of 2-forms $\xi, \eta$ is defined by $\langle 
\xi , \eta \rangle := \frac{1}{2} \xi (X_k , X_l ) \cdot \eta (X^k , X^l)$.
By (\ref{gl-42}), the trace of (\ref{gl-49}) yields the well known equation
\begin{equation} \label{gl-51}
S_{\star} - S = | \nabla \Omega |^2 . 
\end{equation}

Let $\phi$ be the $(0,2)$-tensor field on $M$ defined by
\begin{displaymath}
\phi (X,Y):= \frac{1}{2} \mathrm{tr} (\nabla_X J \cdot \nabla_{JY} J)=
- \frac{1}{2} \langle \nabla_X J, \nabla_{JY} J \rangle .
\end{displaymath}

Using (\ref{gl-5}) and (\ref{gl-9}) we see that $\phi$ has the properties
\begin{equation} \label{gl-52}
\phi (X,Y)= \phi (JX, JY) = - \phi (Y,X) . 
\end{equation}

Thus, $\phi$ is a $J$-invariant 2-form. Moreover, by definition, it
holds that
\begin{equation} \label{gl-53}
\phi (X, JY)= \langle \nabla_X \Omega , \nabla_Y \Omega \rangle . 
\end{equation}

Gray's identity (\ref{gl-39}) provides
\begin{equation} \label{gl-54}
| \tilde{R}^+ |^2 = \frac{1}{2} | \phi |^2 . 
\end{equation}

By (\ref{gl-33}) and (\ref{gl-34}), this implies
\begin{equation} \label{gl-55}
| \tilde{R} |^2 = \frac{1}{2} | \phi |^2 + | \tilde{R}^-|^2 . 
\end{equation}

An $\Omega$-contraction of equation (\ref{gl-47}) yields
\begin{equation} \label{gl-56}
g(( \nabla^2_{X,X_k} J)JX^k , Y)= \frac{1}{2} \mathrm{tr}
(J \circ \nabla_{X,Y} J) . 
\end{equation}

Furthermore, by (14) we find
\begin{equation} \label{gl-57}
\frac{1}{2} \mathrm{tr} (J \circ \nabla^2_{X,Y} J)= \phi (X, JY) . 
\end{equation}

On the other hand, using (\ref{gl-14}), (\ref{gl-15}), (\ref{gl-16})
we have 
\begin{equation} \label{gl-58}
(\nabla^2_{X,X_k} J)JX^k = - (\nabla_{X_k} J \circ \nabla_X J) X^k . 
\end{equation}

By (\ref{gl-56}) and (\ref{gl-57}), this provides the identity
\begin{equation} \label{gl-59}
g(( \nabla_{X_k} J \circ \nabla_X J) X^k , Y)= - \phi (X, JY) . 
\end{equation}

\section{\bf Weitzenb\"ock formulas}

Let $(M, g,J)$ be any almost K\"ahler manifold. By Proposition 1 in 
\cite{5} and (\ref{gl-55}), there is a vector field $V_1$ on $M$ such
that
\begin{equation} \label{gl-60}
\frac{1}{2} | \nabla^* \nabla \Omega |^2 + | \tilde{R} |^2 - | \Ric^- |^2 +2
\langle \rho, \nabla^* \nabla \Omega \rangle - 2 \langle \rho , \phi \rangle +
\mathrm{div} (V_1) =0.
\end{equation}

In contrast to the paper \cite{5} the definition $| \tilde{R} |^2 :=
\sum\limits_{k,l} | \tilde{R} (X_k , X_l )|^2$ is used here. The authors
show that Sekigawa's theorem is an immediate consequence of the basic
Weitzenb\"ock formula (\ref{gl-60}). In the following we modify this
formula in order to obtain some more general results.\\

\begin{lemma} \label{lem-2-1}
{\it For any almost K\"ahler manifold, we have the equations
\begin{equation} \label{gl-61}
2 \langle \rho , \phi \rangle = \langle \nabla_{\Ric (X_k)} \Omega, 
\nabla_{X^k} \Omega \rangle , 
\end{equation}
\begin{equation} \label{gl-62}
\langle \rho , \nabla^* \nabla \Omega \rangle = 2 \langle \Ric , 
\widetilde{\Ric} \rangle . 
\end{equation}}
\end{lemma}

\vspace{0.3cm}

\begin{proof}
We calculate \\

$\displaystyle
2 \langle \rho, \phi \rangle = \rho (X^k , X^l) \cdot \phi (X_k , X_l)= 
\frac{1}{2} g ((J \circ \Ric^+)X^k , X^l) \cdot \mathrm{tr} (\nabla_{X_k}
J \circ \nabla_{J X_l}J )=$\\

$ \displaystyle
=\frac{1}{2} \mathrm{tr} (\nabla_{X_k} J \circ \nabla_{J ((J \circ \Ric^+
) X^k)} J)= \phi (X_k , (J \circ \Ric^+ )X^k) \stackrel{(\ref{gl-52})}{=}$\\

$ \displaystyle
= \phi (X_k , (J \circ \Ric ) X^k) \stackrel{(\ref{gl-53})}{=} \langle
\nabla_{X_k} \Omega , \nabla_{\Ric (X^k)} \Omega \rangle . $\\

This proves (\ref{gl-61}). Further, we have
\begin{displaymath}
\langle \rho , \nabla^* \nabla \Omega \rangle = \frac{1}{2} \langle
J \circ \Ric^+ , 
\nabla^* \nabla J \rangle = - \frac{1}{4} \langle \Ric^+ , \{ \nabla^* \nabla
J , J \} \rangle \stackrel{(\ref{gl-45})}{=} 2 \langle \Ric^+ , 
\widetilde{\Ric} \rangle . 
\end{displaymath}

This yields (\ref{gl-62}) since $\langle \Ric^- , \widetilde{\Ric} \rangle
=0$ by (\ref{gl-24}) , (\ref{gl-44}).
\end{proof}

\vspace{0.3cm}

We introduce the vector-valued 2-form $\psi$ defined by
\begin{displaymath}
\psi (X,Y) := \frac{1}{8} ([ \nabla_X \Ric , J]Y - [\nabla_Y \Ric , J] X-
[\nabla_{JX} \Ric , J] JY + [\nabla_{JY} \Ric , J] JX ) . 
\end{displaymath}

A simple calculation using (\ref{gl-36}) , (\ref{gl-37}) provides the
equation
\begin{equation} \label{gl-63}
\frac{1}{4} g(J( \varphi (X^k , X^l)), (\nabla_{X_k} \Ric ) X_l - (
\nabla_{X_l} \Ric ) X_k)= \langle \varphi , \psi \rangle
\end{equation}

with $\langle \varphi , \psi \rangle := \frac{1}{2} g( \varphi (X^k , X^l ) ,
\psi (X_k , X_l )) . $\\

\begin{lemma} \label{lem-2-2}
{\it For any almost K\"ahler manifold, it holds that 
\begin{equation} \label{gl-64}
2 \langle \Ric , \widetilde{\Ric} \rangle = 2 \langle \rho , \phi \rangle +
2 \langle \varphi , \psi \rangle + | \Ric^- |^2 + \mathrm{div} (V_2) , 
\end{equation}

where $V_2$ is the vector field locally defined by
\begin{displaymath}
V_2 := g (\Ric (X^l) , ( J  \circ \nabla_{X_l} J) X^k ) \cdot X_k . 
\end{displaymath}}
\end{lemma}

\vspace{0.3cm}

\begin{proof}
We calculate\\

$\displaystyle
2 \langle \Ric , \widetilde{\Ric} \rangle \stackrel{(\ref{gl-49})}{=}
\langle \Ric , \Ric^+_{\star} - \Ric^+ \rangle \stackrel{(\ref{gl-50})}{=}$\\

$\displaystyle
= \langle \Ric, \Ric^-_{\star} + \Ric^- \rangle + g( \Ric (X^l) , (J \circ
\nabla_{X_k, X_l} J) X^k) $.\\

By $\langle \Ric , \Ric^-_{\star} \rangle = \langle \Ric^+ , \Ric^- \rangle
=0$, this yields\\

$(*) \hfill 2 \langle \Ric , \widetilde{\Ric} \rangle = | \Ric^- |^2 + 
g(\Ric (X^l) , (J \circ \nabla_{X_k , X_l} J)X^k)$ \hfill \mbox{}\\

Let $x \in M$ be any point and $(X_1 , \ldots , X_n)$ an orthonormal
frame in neighbourhood of $x$ with $(\nabla X_k)_x =0$ for $k= 1 , \ldots ,n$.
Then, at $x \in M$, we have
\begin{eqnarray*}
&& g (\Ric (X^l) , (J \circ \nabla^2_{X_k , X_l} J)X^k)= X_k (g(\Ric
(X^l), (J \circ \nabla_{X_l} J) X^k)) - \\
&& -g (\Ric (X^l) , (\nabla_{X_k} J \circ \nabla_{X_l} J) X^k) - g((
\nabla_{X_k} \Ric) X^l , (J \circ \nabla_{X_l} J) X^k)
\stackrel{(\ref{gl-59}), (5), (9)}{=}\\
&=& \mathrm{div} (V_2) + \phi (X_l , (J \circ \Ric )X^l) - g (J((\nabla_{X^l}
J)X^k - (\nabla_{X^k} J) X^l) , ( \nabla_{X_k} \Ric )X_l )
\stackrel{(\ref{gl-53})}{=}\\
&=& \mathrm{div} (V_2) + \langle \nabla_{X_l} \Omega , \nabla_{\Ric (X^l)}
\Omega \rangle + g(J(\varphi (X^k , X^l)), (\nabla_{X_k} \Ric) X_l)
\stackrel{(\ref{gl-61})}{=}\\
&=& \mathrm{div} (V_2) + 2 \langle \rho , \phi \rangle + \frac{1}{2}
g(J(\varphi (X^k , X^l)), (\nabla_{X_k} \Ric )X_l - (\nabla_{X_l} \Ric)
X_k) \stackrel{(\ref{gl-63})}{=}\\
&=& \mathrm{div} (V_2) + 2 \langle \rho , \phi \rangle + 2 \langle 
\varphi , \psi \rangle
\end{eqnarray*}

and, hence,\\

$(2*) \hfill g(\Ric (X^l) , (J \circ \nabla^2_{X_k , X_l} J) X^k)=
2 \langle \rho , \phi \rangle + 2 \langle \varphi , \psi \rangle +
\mathrm{div} (V_2)$. \hfill \mbox{}\\

Inserting $(2*)$ into $(*)$ we obtain (\ref{gl-64}).
\end{proof}

\vspace{0.3cm}

Using (\ref{gl-48}), (\ref{gl-49}) and $\langle \Ric , \widetilde{\Ric} 
\rangle = \langle \Ric^+ , \widetilde{\Ric} \rangle$ we obtain by 
straightforward calculations\\

\begin{lemma} \label{lem-2-3}
{\it The identities
\begin{equation} \label{gl-65}
\frac{1}{2} | \nabla^* \nabla \Omega |^2 = | \Ric^-_{\star} |^2 + 4|
\widetilde{\Ric} |^2 , 
\end{equation}
\begin{equation} \label{gl-66}
| \widetilde{\Ric} |^2 + \langle \Ric , \widetilde{\Ric} \rangle =
\frac{1}{4} (| \Ric^+_{\star} |^2 - | \Ric^+ |^2 ) , 
\end{equation}
\begin{equation} \label{gl-67}
2 | \widetilde{\Ric} |^2 + \langle \Ric , \widetilde{\Ric} \rangle =
\langle \Ric^+_{\star} , \widetilde{\Ric} \rangle
\end{equation}

are valid for any almost K\"ahler manifold.}
\end{lemma}

\vspace{0.3cm}

\begin{lemma} \label{lem-2-4}
{\it Let $(M,g,J)$ be any almost K\"ahler manifold. Then we have the equations
\begin{equation} \label{gl-68}
| \tilde{R} |^2 + | \Ric_{\star} |^2 - | \Ric |^2 - 2 \langle \rho , \phi
\rangle + \mathrm{div} (V_1) =0 , 
\end{equation}
\begin{equation} \label{gl-69}
| \tilde{R} |^2 + | \Ric^-_{\star} |^2 + 2 \langle \Ric^+_{\star} ,
\widetilde{\Ric} \rangle + 2 \langle \varphi , \psi \rangle +
\mathrm{div} (V_1 + V_2)=0 . 
\end{equation}}
\end{lemma}

\vspace{0.3cm}

\begin{proof}
Inserting (\ref{gl-62}) and (\ref{gl-65}) into (\ref{gl-60}) we find
(\ref{gl-68}) using (\ref{gl-66}). We eliminate the term $2 \langle \rho , 
\phi \rangle $ in (\ref{gl-68}) by (\ref{gl-64}) and  then we obtain 
(\ref{gl-69}) using (\ref{gl-66}) , (\ref{gl-67}).
\end{proof}

\vspace{0.3cm}

Let $(M,g,J)$ be any almost Hermitian $n$-manifold. Then we consider the 
function $q (J)$ locally defined by
\begin{displaymath}
q(J) := g(( \nabla^2_{X_k , X_l} \Ric ) J X^k , JX^l ).
\end{displaymath}

Obviously, it holds that
\begin{equation} \label{gl-70}
q(J) = \frac{1}{2} g(( \nabla^2_{X_k , X_l} \Ric + \nabla^2_{X_l , X_k} \Ric)
JX^k , JX^l)= \frac{1}{2} g(( \nabla^2_{JX_k , JX_l} \Ric + 
\nabla^2_{JX_l , JX_k} \Ric) X^k , X^l ) . 
\end{equation}

If $M$ is compact, then we consider the number
\begin{displaymath}
Q(J) := \int_M q (J) \omega , 
\end{displaymath}

where $\omega := \frac{1}{m!} \Omega^m (n=2m)$ is the volume form. Since
in the K\"ahler case $[\nabla^2 \Ric , J]=0$ , i.e., $[\nabla^2_{X,Y}
\Ric , J]=0$ for all vector fields $X$ and $Y$, we have
\begin{displaymath}
q(J) = g(( \nabla^2_{X_k , X_l} \Ric ) X^k , X^l)= - \frac{1}{2} \Delta S
\end{displaymath}

then and, hence, $Q(J)=0$. Thus, $Q(J)$ is an obstruction against the 
K\"ahler property of any compact almost Hermitian manifold. Furthermore,
if $(M, g,J)$ is almost K\"ahler, then a simple calculation yields
\begin{equation} \label{gl-71}
q (J) = 2 \langle \varphi , \psi \rangle + \mathrm{div} (V_3) , 
\end{equation}

where $V_3$ is the vector field on $M$ locally given by $V_3 := g ((
\nabla_{X_l} \Ric ) JX^l , JX^k) \cdot X_k$. This implies
\begin{equation} \label{gl-72}
Q(J) = 2 \int_M \langle \varphi , \psi \rangle \omega . 
\end{equation}

By (\ref{gl-72}), Lemma \ref{lem-2-4} provides immediately\\

\begin{proposition} \label{prop-2-1}
{\it For any compact almost K\"ahler manifold, the equations
\begin{equation} \label{gl-73}
\int_M (| \tilde{R}|^2 + | \Ric_{\star} |^2 - | \Ric |^2 - 2 \langle \rho , 
\phi \rangle ) \omega =0 , 
\end{equation}
\begin{equation} \label{gl-74}
Q(J) + \int_M (| \tilde{R} |^2 + | \Ric^-_{\star} |^2 + 2 \langle 
\Ric^+_{\star} , \widetilde{\Ric} \rangle ) \omega =0
\end{equation}

are valid.\\}
\end{proposition}

\section{\bf Applications}

For any almost K\"ahler manifold, it is very natural to compare the tensors
$\Ric^+_{\star}$ and $\Ric^+$ since both tensors are symmetric, Hermitian
and coincide with the Ricci tensor in the K\"ahler case. In particular,
a necessary condition for an almost K\"ahler manifold to be K\"ahler is
that $\Ric^+_{\star}$ and $\Ric^+$ have the same length $(|\Ric^+_{\star}
|= | \Ric^+ |)$.\\

\begin{proposition} \label{prop-3-1}
{\it  Let $(M,g,J)$ be an almost K\"ahler manifold with $| \Ric^+_{\star}| =
| \Ric^+|$. Then $(M,g,J)$ is K\"ahler if $\Ric$ or $\Ric^+$ is
semi-positive or if $\Ric^+_{\star}$ is semi-negative.}
\end{proposition}

\vspace{0.3cm}

\begin{proof}
Since $\langle \Ric , \widetilde{\Ric} \rangle = \langle \Ric^+ , 
\widetilde{\Ric} \rangle$ and $\widetilde{\Ric} \ge 0$, our supposition
$\Ric \ge 0$ or $\Ric^+ \ge 0$, respectively, implies $\langle \Ric , 
\widetilde{\Ric} \rangle \ge 0$. Thus, the assumption $|\Ric^+_{\star} |=
| \Ric^+ |$ forces $\widetilde{\Ric} =0$ by (66). Using (66) and (67)
we immediately obtain the equation
\begin{equation} \label{gl-75}
\langle \Ric^+_{\star} , \widetilde{\Ric} \rangle - | \widetilde{\Ric}
|^2 = \frac{1}{4} (| \Ric^+_{\star} |^2 - | \Ric^+ |^2 ) . 
\end{equation}

This shows that the suppositions $| \Ric^+_{\star} |= | \Ric^+ |$ and 
$\Ric^+_{\star} \le 0$ also force $\widetilde{\Ric} =0$. By (41), we
have $\mathrm{tr} (\widetilde{\Ric})= \frac{1}{4} | \nabla J |^2$. Thus,
$\widetilde{\Ric} =0$ implies $\nabla J=0$.
\end{proof}

\vspace{0.3cm}

In dimension 4, we have the following stronger result.\\

\begin{proposition} \label{prop-3-2}
{\it Let $(M,g,J)$ be an almost K\"ahler 4-manifold such that $\Ric^+_{\star}$
and $\Ric^+$ have the same length. Then $(M,g,J)$ is K\"ahler if $S \ge
0$ or $S_{\star} \le 0$ or if $\mathrm{supp} (S + S_{\star})=M$.}
\end{proposition}

\vspace{0.3cm}

\begin{proof}
It is known that, for any almost K\"ahler 4-manifold, the endomorphism
$\nabla_{X_k} J \circ \nabla_{X^k} J$ is a multiple of the identity
map \cite{28}. By (41), this yields
\begin{equation} \label{gl-76}
\widetilde{\Ric} = \frac{1}{8} | \nabla \Omega |^2 . 
\end{equation}

Inserting (76) into (66) and (75) we obtain the equations
\begin{equation} \label{77}
4(| \Ric^+_{\star} |^2 - | \Ric^+ |^2 )= | \nabla \Omega |^4 + 2
S \cdot | \nabla \Omega |^2 = 2 S_{\star} \cdot | \nabla \Omega |^2 - 
| \nabla \Omega |^4 , 
\end{equation}

which provide
\begin{equation} \label{gl-78}
4 (| \Ric^+_{\star} |^2 - | \Ric^+ |^2 )=(S + S_{\star} ) \cdot | 
\nabla \Omega |^2 . 
\end{equation}

Obviously, our suppositions imply $\nabla \Omega =0$ by (77), (78).
\end{proof}

\vspace{0.3cm}

Another necessary condition for an almost K\"ahler manifold to be K\"ahler is
$| \Ric_{\star} |= | \Ric|$. The following theorem shows that this
condition is also sufficient in the compact case with semi-negative Ricci 
tensor $(\Ric \le 0)$.\\

\begin{theorem} \label{thm-3-1}
{\it Let $(M,g,J)$ be any compact almost K\"ahler manifold such that at least
one of the tensors $\Ric, \Ric^+$ or $\Ric^+_{\star}$ is semi-negative.
Then $(M,g,J)$ is K\"ahler if $\Ric_{\star}$ and $\Ric$ have the same length.}
\end{theorem}

\vspace{0.3cm}

\begin{proof}
It holds that
\begin{eqnarray*}
\langle \nabla_{\Ric (X_k)} \Omega , \nabla_{X^k} \Omega \rangle 
\stackrel{(5),(9)}{=} \langle \nabla_{\Ric^+ (X_k)} \Omega , 
\nabla_{X^k} \Omega \rangle \stackrel{(\ref{gl-49})}{=}
 \langle \nabla_{\Ric^+_{\star} (X_k)} \Omega , \nabla_{X^k}
\Omega \rangle - 2 \langle \nabla_{\widetilde{\Ric} (X_k)} \Omega , 
\nabla_{X^k} \Omega \rangle . 
\end{eqnarray*}

By (\ref{gl-61}), this yields
\begin{equation} \label{gl-79}
2 \langle \rho , \phi \rangle = \langle \nabla_{\Ric^+ (X_k)} \Omega , 
\nabla_{X^k} \Omega \rangle = \langle \nabla_{\Ric^+_{\star} (X_k)}
\Omega , \nabla_{X_k} \Omega \rangle - 2 \langle \nabla_{\widetilde{\Ric}
(X_k)} \Omega , \nabla_{X^k} \Omega \rangle . 
\end{equation}

Using (\ref{gl-61}) and (\ref{gl-79}) we see that each of the 
suppositions $\Ric \le 0, \Ric^+ \le 0$ or $\Ric^+_{\star} \le 0$ implies
$\langle \rho , \phi \rangle \le 0$. Thus, by (\ref{gl-55}) and 
(\ref{gl-73}), we obtain the inequality
\begin{equation} \label{gl-80}
\int_M (\frac{1}{2} | \phi |^2 + | \tilde{R}^- |^2 + | \Ric_{\star} |^2
- | \Ric |^2 ) \omega \le 0 , 
\end{equation}

which provides $\phi =0$ if $| \Ric_{\star} |= | \Ric |$. By (\ref{gl-53}), 
$\phi =0$ implies $\nabla \Omega =0$.
\end{proof}

\vspace{0.3cm}

The application of the integral formula (\ref{gl-73}) to the Einstein case
yields\\

\begin{theorem} \label{thm-3-2}
{\it A compact Einstein almost K\"ahler $n$-manifold is necessarily K\"ahler if
the inequality
\begin{equation} \label{gl-81}
| \tilde{R}^- |^2 + | \Ric_{\star} |^2 \ge \frac{1}{n} S \cdot S_{\star}
\end{equation}

is satisfied.}
\end{theorem}

\vspace{0.3cm}

\begin{proof}
In the Einstein case of $\Ric = \frac{S}{n}$, we have
\begin{displaymath}
| \Ric |^2 + 2 \langle \rho , \phi \rangle \stackrel{(\ref{gl-61})}{=}
\frac{S^2}{n} + \frac{S}{n} |\nabla \Omega |^2 
\stackrel{(\ref{gl-51})}{=} \frac{1}{n} S \cdot S_{\star} . 
\end{displaymath}

Inserting this into (\ref{gl-73}) we see that (\ref{gl-81}) forces
$\nabla \Omega =0$.
\end{proof}

\vspace{0.3cm}

Since inequality (\ref{gl-81}) is equivalent to
\begin{equation} \label{gl-82}
| \tilde{R}^- |^2 + | \Ric_{\star} - \frac{S_{\star}}{n} |^2 +
\frac{S_{\star}}{n} | \nabla \Omega |^2 \ge 0 , 
\end{equation}

we immediately obtain\\

\begin{corollary} \label{cor-3-1}
{\it Every compact Einstein almost K\"ahler manifold with non-negative star
scalar curvature is K\"ahler.}
\end{corollary}

\vspace{0.3cm}

This corollary is a slight generalization of Sekigawa's theorem since
$S \ge 0$ implies $S_{\star} \ge 0$ by (\ref{gl-51}).\\
Concerning the Goldberg conjecture it may be interesting to investigate
inequality (\ref{gl-81}) in suitable geometrical situations with $S<0$.\\
Now we consider the integral formula (\ref{gl-74}). The main results of this
paper are the following two theorems.\\

\begin{theorem} \label{thm-3-3}
{\it Let $(M, g,J)$ be any compact almost K\"ahler manifold such that at least
one of the tensors $\Ric, \Ric^+$ or $\Ric^+_{\star}$ is semi-positive.
Then $J$ is integrable if $Q(J)=0$.}
\end{theorem}

\vspace{0.3cm}

\begin{proof}
By (\ref{gl-43}) and (\ref{gl-67}), each of the suppositions $\Ric \ge 0$, 
$\Ric^+ \ge 0$ or $\Ric^+_{\star} \ge 0$, respectively, implies $\langle
\Ric^+_{\star} , \widetilde{\Ric} \rangle \ge 0$. Thus, (\ref{gl-74}) and
$Q(J)=0$ force $\tilde{R} =0$, and, hence, $\nabla \Omega =0$ by (\ref{gl-55})
and (\ref{gl-53}).
\end{proof}

\vspace{0.3cm}

By (\ref{gl-76}), equation (\ref{gl-74}) yields immediately\\

\begin{theorem} \label{thm-3-4}
{\it If $(M,g,J)$ is a compact almost K\"ahler 4-manifold with
$S_{\star} \ge 0$, then $Q(J)=0$ implies $\nabla J=0$.}
\end{theorem}

\vspace{0.3cm}

\begin{remark} \label{remark-3-1}
The supposition in Theorem \ref{thm-3-3} that at least one of the tensors 
$\Ric, \Ric^+$ or $\Ric^+_{\star}$ is semi-positive can be replaced
by the weaker condition
\begin{equation} \label{gl-83}
| \tilde{R}^- |^2 + | \Ric^-_{\star} |^2 + 2 \langle \Ric^+_{\star} , 
\widetilde{\Ric} \rangle \ge 0 . 
\end{equation}

This condition does not involve any derivatives of the almost complex
structure $J$. By (\ref{gl-49}) , (\ref{gl-83}) becomes
\begin{equation} \label{gl-84}
| \tilde{R}^- |^2 + | \Ric_{\star} |^2 \ge \langle \Ric^+_{\star} , 
\Ric^+ \rangle . 
\end{equation}

Moreover, in dimension 4, (\ref{gl-83}) is equivalent to
\begin{equation} \label{gl-85}
| \tilde{R}^- |^2 + | \Ric^-_{\star} |^2 + \frac{S_{\star}}{4}
| \nabla \Omega |^2 \ge 0 . 
\end{equation}
\end{remark}

In the K\"ahler case, these inequalities are satisfied trivially.\\

\begin{remark} \label{remark-3-2}
(i) The condition
\begin{equation} \label{gl-86}
[ \nabla \Ric , J] =0 , 
\end{equation}

i.e., $[\nabla_X \Ric , J]=0$ for all vector fields $X$, implies
$\psi =0$ and, hence, $Q(J)=0$ by (\ref{gl-72}). Obviously, any K\"ahler
manifold satisfies this condition.\\

(ii) By (\ref{gl-70}), the condition
\begin{equation} \label{gl-87}
[\nabla^2_{X,Y} \Ric + \nabla^2_{Y,X} \Ric , J] =0
\end{equation}

yields $q(J)= - \frac{1}{2} \Delta S$ and, hence, $Q(J) =0$. In 
particular, the stronger supposition that $[\nabla^2 \Ric , J]=0$ forces
$Q(J)=0$. Again, any K\"ahler manifold fulfils  this condition.\\

(iii) The second one of the equations (\ref{gl-70}) shows that
\begin{equation} \label{gl-88}
\nabla^2_{JX, JY} \Ric + \nabla^2_{JY, JX} \Ric = \nabla^2_{X,Y} \Ric
+ \nabla^2_{Y,X} \Ric
\end{equation}

implies $Q(J)=0$.\\

(iv) Using (\ref{gl-63}) a simple calculation yields that the supposition
\begin{equation} \label{gl-89}
(\nabla_X \Ric)Y - (\nabla_Y \Ric )X= \frac{1}{2(n-1)} (X(S)Y-Y(S)X)
\end{equation}

provides $\langle \varphi , \psi \rangle =0$ and, hence, $Q(J)=0$ by
(\ref{gl-72}).
\end{remark}

It is well known that, for dimension $n \ge 4$, (\ref{gl-89}) is
equivalent to the condition that the Weyl tensor $W$ is divergence-free
or co-closed $(\delta W=0)$. Since $\delta W =0$ implies $dW=0$ (second
Bianchi identity for $W$), $W$ is also called harmonic in this case. Examples
of Riemannian manifolds with harmonic Weyl tensor are given in \cite{7},
Chapter 16.D. Such manifolds are also called nearly conformally symmetric.\\

(v) A straightforward calculation shows that (\ref{gl-89}) implies
\begin{equation} \label{gl-90}
\nabla^2_{X,Y} \Ric + \nabla^2_{Y,X} \Ric =( \nabla^2_{X , \cdot }\  \Ric)Y+
(\nabla^2_{Y , \cdot} \ \Ric)X+ \frac{1}{2(n-1)} (2 \nabla^2_{X,Y} S - 
\nabla_X dS \otimes Y - \nabla_Y dS \otimes X) . 
\end{equation}

Inserting this into (\ref{gl-70}) we find $q(J) = - \frac{1}{2(n-1)} 
\Delta S$. Thus, (\ref{gl-90}) forces $Q(J)=0.$\\[0.4cm]

Remark \ref{remark-3-2} yields some possibilities for concrete applications
of our main theorems. The simplest case of application is the situation with
parallel Ricci tensor in which each of the conditions (\ref{gl-86}) - 
(\ref{gl-90}) is satisfied trivially. By Remark \ref{remark-3-2}, (i), 
the following corollaries are immediate consequences of Theorem 
\ref{thm-3-3} and Theorem \ref{thm-3-4}, respectively.\\

\begin{corollary} \label{cor-3-2}
{\it Let $(M,g,J)$ be any compact almost K\"ahler manifold 
such that at least one of
the tensors $\Ric, \Ric^+$ or $\Ric^+_{\star}$ is semi-positive. Then
$[\nabla \Ric , J]=0$ implies $\nabla J=0$.}
\end{corollary}

\vspace{0.3cm}

\begin{corollary} \label{cor-3-3}
{\it Every compact almost K\"ahler 4-manifold $(M, g,J)$ with $S_{\star} \ge 0$
and $[\nabla \Ric , J] =0$ is K\"ahler.}
\end{corollary}

\vspace{0.3cm}

Our Corollary \ref{cor-3-2} is a generalization of Sekigawa's theorem. It
shows that the assertion of this theorem is true not only for compact 
Einstein manifolds of non-negative scalar curvature but also for 
Riemannian products of such manifolds. It seems that condition 
(\ref{gl-89}) is a suitable generalization of
the supposition that 
$\nabla \Ric =0$ for our purposes. But this is not the case 
by the fact that the Ricci tensor of a K\"ahler manifold with harmonic Weyl 
tensor is parallel (\cite{7}, 16.30 Prop.). Thus, every almost K\"ahler
manifold with harmonic Weyl tensor and non-parallel Ricci tensor cannot
be K\"ahler and, hence, is  strictly almost K\"ahler by definition. Taking
into account this fact, by Remark \ref{remark-3-2}, (iv) and  Theorem
\ref{thm-3-3}, we immediately obtain\\

\begin{corollary} \label{cor-3-4}
{\it A compact almost Hermitian manifold of dimension $n \ge 4$ with 
harmonic Weyl tensor and non-parallel, semi-positive Ricci tensor cannot
be almost K\"ahler.}
\end{corollary}

\vspace{0.3cm}

Finally, Theorem \ref{thm-3-4} implies\\

\begin{corollary} \label{cor-3-5}
{\it Every compact almost Hermitian 4-manifold with harmonic Weyl tensor,
non-parallel Ricci tensor and non-negative scalar curvature is not almost
K\"ahler.}
\end{corollary}

\vspace{0.3cm}

In particular, Corollary \ref{cor-3-5} shows that there are no compact 
almost K\"ahler
4-manifolds with harmonic Weyl tensor and non-constant, non-negative
scalar curvature.\\

\def\refname{References}

\end{document}